\documentclass{amsart}
\usepackage{amssymb,hyperref} 
\newtheorem{thm}{Theorem}[section]

\theoremstyle{definition}

\numberwithin{equation}{section}

\begin{document}

\title[Commuting Probability]{Commuting Probability of Compact Groups}%

\author[A. Abdollahi]{Alireza Abdollahi}%
\address{Department of Pure Mathematics, Faculty of Mathematics and Statistics, University of Isfahan, Isfahan 81746-73441, Iran.} 
\email{a.abdollahi@math.ui.ac.ir}%
\author[M. Soleimani Malekan]{Meisam Soleimani Malekan}%
\address{Department of Pure Mathematics, Faculty of Mathematics and Statistics, University of Isfahan, Isfahan 81746-73441, Iran; Institute for Research in Fundamental Sciences, School of Mathematics, Tehran, Iran.} 
\email{msmalekan@gmail.com}


\begin{abstract}
For any (Hausdorff) compact group $G$ with the normalized Haar measure ${\mathbf m}_G$, denote by ${\rm cp}(G)$ the probability ${\mathbf m}_{G\times G}(\{(x,y)\in G\times G \;|\; xy=yx\})$ of commuting a randomly chosen pair of elements of $G$. Here we prove that if ${\rm cp}(G)>0$, then there exists a finite group $H$ such that ${\rm cp}(G)= \frac{{\rm cp}(H)}{|G:F|^2}$, where $F$ is the FC-center of $G$ i.e. the set of all elements of $G$ whose conjugacy classes are finite and $H$ is isoclinic to $F$ with ${\rm cp}(F)={\rm cp}(H)$.  The latter equality  enables one to transfer many existing results concerning commuting probability of finite groups to one of compact groups. For example, here for a compact group $G$ we prove that  if ${\rm cp}(G)>\frac{3}{40}$ then either $G$ is solvable or, else  $G\cong A_5 \times T$ for some abelian group $T$, in which case ${\rm cp}(G)=\frac{1}{12}$; where $A_5$ denotes the alternating group of degree $5$. 
\end{abstract}
\maketitle
\section{\bf Introduction and Results}

Let $G$ be a (Hausdorff) compact group. Then $G$ has a unique normalized Haar measure denoted by ${\mathbf m}_G$. The set ${\mathcal A}(G):=\{(x,y)\in G\times G \;|\; xy=yx\}$ is a measurable set (actually a closed set) of $G\times G$ and its measure in $G\times G$ will be denoted by ${\rm cp}(G)$. In the case $G$ is finite, ${\rm cp}(G)$ has been extensively studied for example see \cite{GR} and \cite{LNY} and references therein. It is famous  (see e.g. \cite{G}) that if $G$ is finite, ${\rm cp}(G)\leq \frac{5}{8}$ whenever $G$ is non-abelian. The latter is generalized for any non-abelian compact group in \cite{G}. The rationality of ${\rm cp}(G)$ is clear in the case of finite $G$ and  for compact groups it is proved in \cite[Theorem 3.11]{HR}. It is also shown in \cite[Theorems 1.1 and 1.2]{HR} that if ${\rm cp}(G)>0$, then $Z(FC(G))$ is open in $G$, where $Z(FC(G))$ is the center of the FC-center $FC(G)$ of $G$, i.e. $FC(G)$ consists of  elements of $G$ whose conjugacy classes are finite. In the latter case the derived subgroup of $FC(G)$ is also finite. We will use these results in the sequel.\\  
Here we prove a formula relating the commuting probability of compact groups to one of finite groups.
\begin{thm}\label{main}
For any compact group $G$ there exists a finite group $H$ such that ${\rm cp}(G)= \frac{{\rm cp}(H)}{|G:FC(G)|^2}$. In particular, if ${\rm cp}(G)>0$, the finite group $H$ can be chosen so that it is isoclinic to $FC(G)$ and ${\rm cp}(FC(G))={\rm cp}(H)$.  
\end{thm} 
Here we are using the usual convention that if $|G:FC(G)|$ is infinite, then $\frac{1}{|G:FC(G)|}=0$.
Recall that, following \cite{H}, two groups $H$ and $K$ are called isoclinic if there exist group isomorphisms $\alpha:H/Z(H)\rightarrow K/Z(K)$ and $\beta: H' \rightarrow K'$ such that for all $x,y\in H$ and all $x'\in (xZ(H))^\alpha$ and $y'\in (yZ(H))^\alpha$, $[x',y']=[x,y]^\beta$, where for a group $T$, $Z(T)$ and $T'$ are denoting the center and the derived subgroups of $T$, respectively. The alternating group of degree $5$ is denoted by $A_5$ as usual.  
In view of Theorem \ref{main} it  is now possible to generalize some of results in finite case for commuting probability to compact groups. For example we give the following.
\begin{thm}\label{main2}
	Let $G$ be a compact group. Then\\
{\rm (1)} \; if ${\rm cp}(G)>\frac{1}{4}$ then both $G'$ and $G/Z(G)$ are finite.\\
{\rm (2)} \; if ${\rm cp}(G)>\frac{3}{40}$ then either $G$ is solvable or, else  $G\cong A_5 \times T$ for some abelian group $T$, in which case ${\rm cp}(G)=\frac{1}{12}$. 
\end{thm}   

The second part of Theorem \ref{main} has already been proved in \cite[Theorem 12]{GR} for finite groups and we will use it in our proof for the compact case. 

\section{\bf Proof of Theorem \ref{main}}

Let $n=|G:Z|$ and $m=|F:Z|$, where $F$ is the FC-center of $G$ and $Z$ is its center.  By \cite[Theorem 1.2]{HR} we may assume that ${\rm cp}(G)>0$ so that $n$ and $m$ are both finite. Suppose that $F=\cup_{i=1}^m x_i Z$ and $G=\cup_{j=1}^n y_j Z$, where $x_i\in F$ and $y_j \in G$.
Let $c_{\ell,k}=1$  if $x_\ell x_k=x_k x_\ell$ and  $c_{\ell,k}=0$ otherwise.\\

We use the following in the sequel.

\noindent (1) \;  If $y\in G\setminus F$, then ${\mathbf m}_G(C_G(y))=0$, where $C_G(y)=\{a\in G \;|\; ay=ya\}$.\\
(2) \; By (1) if $y_j \in G\setminus F$, ${\mathbf m}_G(x_iZ \cap C_G(y_j))=0$. If $y_j \in F$, we may assume without loss of generality that $y_j=x_\ell$ for some $\ell$. \\ 
(3) \;  $x_iZ \cap C_G(x_\ell)=\begin{cases} x_iZ & {\rm if} \; c_{\ell,i}=1 \\ \varnothing & {\rm otherwise} \end{cases}$. \\
(4) For any $x\in F$
$C_{G}(x) \cap y_jZ=\begin{cases} y_jZ & {\rm if} \; x\in C_G(y_j) \\ \varnothing & {\rm otherwise} \end{cases}$.\\

Let ${\mathbf 1}_{{\mathcal A}(G)}$ be the characteristic function of ${\mathcal A}(G)=\{(x,y)\in G\times G \;|\; xy=yx\}$. \\
 The following two first equlities are from \cite{G}. 

\begin{align*}
{\rm cp}(G)&=\int_{G\times G} {\mathbf 1}_{{\mathcal A}(G)}(x,y) d(x,y)=\int_G  {\mathbf m}_G(C_G(x)) dx \;\;\; {\rm by \; Fubini \; theorem} \\ 
&\overset{(1)}{=}\int_F {\mathbf m}_G(C_G(x)) dx  \\
&=\sum_{i=1}^m\int_{x_i Z}  {\mathbf m}_G(C_G(x)) dx = \sum_{i=1}^m \sum_{j=1}^n \int_{x_iZ}  {\mathbf m}_G(C_{G}(x) \cap y_jZ) dx \\
&\overset{(2)}{=}\sum_{i=1}^m \sum_{j=1}^n \int_{x_iZ \cap C_G(y_j)}  {\mathbf m}_G(y_jZ) dx =\sum_{i=1}^m \sum_{j=1}^n \int_{x_iZ \cap C_G(y_j)}  {\mathbf m}_G(Z) dx\\
&={\mathbf m}_G(Z) \sum_{i=1}^m \sum_{j=1}^n \int_{x_iZ \cap C_G(y_j)} {\mathbf 1}_{x_iZ \cap C_G(y_j)}(x)  dx ={\mathbf m}_G(Z) \sum_{i=1}^m \sum_{j=1}^n {\mathbf m}_G(x_iZ \cap C_G(y_j))   \\
&\overset{(2)}{=}{\mathbf m}_G(Z) \sum_{i,j=1}^m {\mathbf m}_G(x_iZ \cap C_G(x_j)) \overset{(3)}{=}{\mathbf m}_G(Z) \sum_{i,j=1}^m {\mathbf m}_G(Z) c_{i,j}   \\
&={\mathbf m}_G(Z)^2 \sum_{i,j=1}^m  c_{i,j} = \frac{1}{|G:Z|^2} \sum_{i,j=1}^m c_{i,j}  
\end{align*}
The latter, in particular, shows the rationality of ${\rm cp}(G)$. 
On the other hand,
\begin{align*}
{\rm cp}(F)&= {\mathbf m}_{F\times F}\big(\bigcup_{i,j=1}^m \{x_iZ \times x_j Z \;|\; x_i x_j = x_j x_i\}\big)\\
 &=\sum_{i,j=1}^m {\mathbf m}_{F\times F}(x_i Z \times x_j Z) c_{i,j} \\
&= \sum_{i,j=1}^m {\mathbf m}_F(Z)^2 c_{i,j} ={\mathbf m}_F(Z)^2 \sum_{i,j=1}^m  c_{i,j}
\end{align*}
It follows that 
\begin{align} 
\sum_{i,j=1}^m  c_{i,j}=\frac{{\rm cp}(F)}{{\mathbf m}_F(Z)^2} \tag{*}
\end{align}
 Thus we also find the following relation.
 $${\rm cp}(G)= \frac{|F:Z|^2}{|G:Z|^2} {\rm cp}(F)=\frac{{\rm cp}(F)}{|G:F|^2}.$$
 By \cite[p. 135, paragraph 4]{H}, there exists a group $H$ such that $Z(H)\leq H'$ and $H$ is isoclinic to $F$. Note that since $F/Z(F)$ and $F'$ are finite, $H$ is a finite group. If $d_{i,j}$ are defined similar to $c_{i,j}$ for $H$, then the isoclinism of $H$ and $F$ implies that  $\sum_{i,j=1}^m  c_{i,j}=\sum_{i,j=1}^m  d_{i,j}$. Note that the same relation as (*) holds for $H$ as well, thus
  $$\frac{{\rm cp}(F)}{{\mathbf m}_F(Z)^2}=\frac{{\rm cp}(H)}{{\mathbf m}_H(Z(H))^2}.$$ 
  Since ${\mathbf m}_F(Z)=|F:Z|$ and ${\mathbf m}_H(Z(H))=|H:Z(H)|$, the isoclinism again implies that 
  ${\rm cp}(F)={\rm cp}(H)$.
This completes the proof of Theorem \ref{main}. 

\section{\bf Proof of Theorem \ref{main2}}

 By Theorem \ref{main} there exists a finite group $H$ isoclinic to $F:=FC(G)$ such that ${\rm cp}(G)=\frac{{\rm cp}(H)}{|G:F|^2}$.\\

(1) \; Suppose that ${\rm cp}(G)>\frac{1}{4}$.  If $G\neq F$,  then $\frac{1}{|G:F|^2}\leq \frac{1}{4}$ and so ${\rm cp}(H)>1$ which is impossible. Thus $G=F$. Now \cite[Theorem 3.11]{HR} follows that both $G'$ and $G/Z(G)$ are finite. \\

(2) \; Suppose that ${\rm cp}(G)>\frac{3}{40}$.  It follows that ${\rm cp}(H)>\frac{3}{40}$ and $|G:F|\in\{1,2,3\}$ so that $G'\leq F$. Now 
\cite[Theorem 12]{GR} implies that $H$ is solvable or $H\cong A_5 \times T$ for some abelian group $T$. If $H$ is solvable, by isoclinism, $H'\cong F'$ and so $F$ is solvable and so is $G$, since $G/F$ is cyclic. Now assume that $H\cong A_5 \times T$ for some abelian group $T$ so that ${\rm cp}(H)=\frac{1}{12}$. Now (1) implies that $G=F$. By isoclinism, $G/Z(G)\cong A_5$ and $G'\cong A_5$. Therefore $G=G'Z(G)$ and so $G\cong A_5 \times Z(G)$.

\end{document}